\documentclass[12pt]{article}
\usepackage{amsmath,amssymb,amsthm,amsfonts}
\usepackage{latexsym}
\begin{document}

\newtheorem{prop}{Proposition}[section]
\newtheorem{cor}{Corollary}[section] 
\newtheorem{theo}{Theorem}[section]
\newtheorem{lem}{Lemma}[section]
\newtheorem{rem}{Remark}[section]
\newtheorem{con}{Conjecture}[section]
\newtheorem{as}{Assumption}[section]
\newtheorem{de}{Definition}[section]

\setcounter{page}{1}
\renewcommand{\theequation}{\thesection.\arabic{equation}}

\begin{center}
{\Large \bf Diffusion approximation for the 
components in critical inhomogeneous random graphs of rank 1.}
\end{center}

\begin{center}    
TATYANA S. TUROVA\footnote{
Mathematical Center, University of
Lund, Box 118, Lund S-221 00, 
Sweden.

\ \ \ Research was carried out at the Mittag-Leffler Institute, Djursholm, Sweden, and partially
supported by the Swedish Natural Science Research
Council.} 
\end{center}

\begin{abstract}
Consider the random graph on $n$ vertices $1, \ldots, n$. Each vertex $i$ is assigned a type $X_i$ with $X_1, \ldots , X_n$ being independent identically distributed as a nonnegative discrete random variable $X$. 
We assume that ${\bf E} X^3<\infty$.
Given types of all vertices, an edge exists between vertices $i$ and $j$ independent of anything else and with probability $\min \{1,
\frac{X_iX_j}{n}\left(1+\frac{a}{n^{1/3}} \right) \}$. We study the critical phase, which is known to take place when ${\bf E} X^2=1$. 
We prove that normalized by $n^{-2/3}$ the asymptotic joint distributions of component sizes of the graph equals the joint distribution of the excursions of a reflecting Brownian motion $B^a(s)$ with diffusion coefficient
$\sqrt{{\bf E}X{\bf E}X^3}$ and drift $a-\frac{{\bf E}X^3}{{\bf E}X}s$. 
This shows that finiteness of ${\bf E}X^3$  is the necessary condition for the diffusion limit.
In particular, we conclude that 
the
size of the largest connected component is of order $n^{2/3}$. 

\end{abstract}


\medskip

\section{Introduction.}
\setcounter{equation}{0}
\subsection{The Model.}
We study here a rank 1 case of a general inhomogeneous
random graph model  introduced  by Bollob{\'a}s,   Janson and   Riordan 
\cite{BJR}. We shall define a
random graph
$G^{\cal V}(n)$ 
with a vertex space
$${\cal V}=(S,\mu, (x_1^{(n)}, \ldots, x_{n}^{(n)})_{n\geq 1}),$$
where $S=\{0,1, \ldots\}$ and $\mu$ is a probability
 on $S$.
No
relationship is assumed between $x_i^{(n)}$ and $x_i^{(n')}$, but to
simplify notations we
shall write further $(x_1, \ldots, x_{n})=(x_1^{(n)}, \ldots, x_{n}^{(n)})$.
For each $n$ let
$(x_1, \ldots, x_{n})$ be
a deterministic or random
sequence  of points in $S$, such that for any
 $A\subseteq S$
\begin{equation}\label{set}
    \frac{\#\{i: x_i\in A\}}{n}\stackrel{P}{\rightarrow}   \mu (A).
\end{equation}
Given the sequence $x_1, \ldots, x_{n}$, we let  $G^{\cal V}(n)$
 be the random graph on $\{1, \ldots, n\}$, such that any
two vertices $i$ and $j$ are connected by an edge independently of the others
and with a probability
\begin{equation}\label{pe}
p_n(x_i, x_j)= \min\left\{ \frac{x_i x_j}{n}\left(1+\frac{a}{n^{1/3}}\right), \ 1 \right\},
\end{equation}
where $a$ is any fixed real constant.

Let $X$ denote a random variable with values in $S$ and probability function $\mu$. It is proved in \cite{BJR} that the phase transition occurs when
${\bf E}X^2=1.$
Our aim here is to derive the asymptotic behaviour of the sizes of the connected components. 

We shall assume that 
\begin{equation}\label{1}
{\bf E}X^2=1,
\end{equation}
and
\begin{equation}\label{2}
{\bf E}X^{3}< \infty.
\end{equation}

Let  $x_1, \ldots , x_n$ be $ i.i.d.$ as random variable $X$.
Then by assumption (\ref{2})
there is an increasing unbounded function $\omega(x)$ such that
bound
\begin{equation}\label{A1}
\max_{1\leq i \leq n}x_i\leq \frac{n^\frac{1}{3}}{\omega(n)},
\end{equation}
holds with probability tending to one as $n\rightarrow \infty$.

\begin{rem} It will be clear that one can extend all the results for the case of non - i.i.d. random variables, assuming, however, 
(\ref{A1}) and
some uniformity in the convergence (\ref{set}).
\end{rem}

Let $C_1(n), C_2(n), \ldots$ denote the ordered sizes of the connected components in $G^{\cal V}(n)$ with $C_1(n)$ being the largest one. We shall find the weak limit of 
\[n^{-2/3}\left( C_1(n), C_2(n), \ldots \right).\]

Let $(W(s), s\geq 0)$ be the standard Brownian motion. Define
\begin{equation}\label{11}
{ W}^a(s)=\sqrt{{\bf E}X{\bf E}X^3}\ W(s)+as-\frac{{\bf E}X^3}{2{\bf E}X}s^2.
\end{equation}
Let $\gamma_1, \gamma_2, \ldots$ denote the ordered lengths of the excursions of the  process
\begin{equation}\label{3}
B(s)={ W}^a(s)- \min_{0\leq s'\leq s}{W}^a(s'), \ \ \ s\geq 0.
\end{equation}

To formulate the convergence result define ${\it l}^2$ to be the set of infinite sequences ${\bf x}=(x_1,  x_2, \ldots)$ with $x_1\geq  x_2 \geq  \ldots\geq 0$
and $\sum_i x_i^2 <\infty$, and give ${\it l}^2$ metric $d({\bf x}, {\bf y})=\sqrt{(x_i-y_i)^2}$.

\begin{theo}\label{C1} The convergence in distribution
\begin{equation}\label{4}
n^{-2/3}\left( C_1(n), C_2(n), \ldots \right)\stackrel{d}{\rightarrow}
\left(\gamma_1, \gamma_2, \ldots\right)
\end{equation}
holds with respect to the ${\it l}^2$ topology.
\end {theo}

The statement of this result is very much inspired by the 
earlier work of  Aldous \cite{A}, who was the first to prove rigorously this theorem not only for the case of homogeneous random graph $G_{n,p}$ (in our setting this corresponds the case when $x_i \equiv 1$), but for  
a particular nonuniform case as well. 
Notice here, that Theorem \ref{C1} and  
the  result  of Proposition 4 in \cite{A} for  a nonuniform random graph
only partially overlap: there are cases which are covered by one theorem, 
but not by another. 
It should be also mentioned that
in the nonuniform case  Aldous derives  in \cite{A} 
asymptotic for the {\it sums of types of vertices} in the components, while  
we state the result directly for the components sizes. However, this is a minor difference, since most likely 
both objects behave in a similar fashion: the critical values
for 
the phase transition 
coincide, and the phase transitions are qualitatively very similar
(at least for the case of $X$ with exponential tail  it is explicitly derived in \cite{T}).

The main difference between
Theorem \ref{C1} and  result of Aldous for a nonuniform random graph
(Proposition 4 in \cite{A}) is in the relations between the assumptions
on the graph model
 and the coefficients of the process $W^a$.
In the proof of Theorem \ref{C1} we 
combine approach of 
 Aldous in \cite{A}  with 
idea of  Martin-L{\"o}f in \cite{M} on a construction of 
 diffusion approximation for a critical epidemics. This allows us to 
derive in a straightforward way the coefficients of the diffusion process $W^a$,
which gives a number of advantages.
In particular, it is transparent in our proof (see Remark \ref{cE}) that
 the  phase transition occurs indeed when  ${\bf E} X^2=1$, and that the famous scaling $n^{2/3}$ is the proper one here.

Theorem \ref{C1} places all the rank 1 graphs with a finite third  moment into the same universality class as a homogeneous $G_{n,p}$ model, as long as the scaling $n^{2/3}$ concerns. This fact was also observed 
by van der Hofstad \cite{H} who (under some additional
assumptions about the distribution of $X$)
obtained good bounds for the probability of order (in $n$)
of the largest  component, and
classified possible critical scalings when only ${\bf E}X^2<\infty$.

Our proof also shows   that finiteness of ${\bf E}X^3$  is the 
necessary condition for the diffusion limit, and it  indicates the 
possibility to find another scaling and a corresponding weak limit when ${\bf E}X^3<\infty$ 
is not the case.

\subsection{The breadth-first walk for inhomogeneous random graph.}
Following ideas from \cite{A} we construct here a process associated with revealing connected components in  inhomogeneous random graph.
The basic procedure is the same as in homogeneous random graphs: 

Given a graph $G^{\cal V}(n)$ and a set 
$$
V_n=\{x_1, \ldots , x_n\}$$
choose a random (size-biased, as we explain later)  vertex in $\{1\ldots , n\}$, and mark it by  $v_1$. Reveal  all
the vertices connected to the marked vertex $v_1$ in the graph $ G^{\cal V}(n)$. 
If the set of 
non-marked revealed vertices 
is not empty pick a random uniform vertex among this set,
mark it $v_2$, and find all the vertices
connected to it but which have not been used previously. 
We continue this process until we end
up with a tree of marked vertices. Then choose again randomly size-biased vertex  among the unrevealed ones and start the same process again until we use all the vertices of the graph. 

We shall now introduce a Markov chain
\begin{equation}\label{MC}
\left({\bar U}(i), {\bar I}(i), x(i)\right), \ \ \ 1\leq i\leq n,
\end{equation}
associated with
 this algorithm, where at any step $i\geq 1$

${\bar U}(i)=(U^x(i), x\in S)$ and $U^x(i)$ denotes the number of unrevealed vertices of type $x$;

${\bar I}(i)=(I^x(i), x\in S)$ and $I^x(i)$ denotes the number of revealed but non-marked  vertices of type $x$;

$x(i)=x_{v_i}$ is the type of the vertex marked at step $i$. 

Notice that all these variables depend on $n$. Revealed but non-marked  vertices at step $i$ we shall simply call
''active vertices at step $i$''.

For the initial state we set
\begin{equation}\label{5}
\begin{array}{ll}
I^x(1)=& 0,\\ \\
U^x(1)=&\#\{1\leq i\leq n: x_i=x\}, \ \ x\in S,
\end{array}
\end{equation}
and let 
$x(1)$ be a random variable with {\it size-biased} distribution 
\begin{equation}\label{defp}
{\bf P}\{x(1)=x\mid {\bar U}(1) \}=\frac{x U^x(1)}{\sum _{x\in S}xU^x(1)}.
\end{equation}

Denote further for any $i\geq 1$
\[I(i)=\sum_{x\in S}I^x(i),\]
\[U(i)=\sum_{x\in S}U^x(i).\]
Then for any $i\geq 1 $ conditionally on $(x(i), {\bar I}(i),{\bar U}(i))$  
the number of  new
''$x$ - neighbours'' of $v_{i}$ in the graph  $G^{\cal V}(n)$ 
is distributed as
\begin{equation}\label{10*}
{\cal N}_n^x(i) \in Bin \Big(U^x(i)-{\bf 1}_{I(i)=0}{\bf 1}_{x(i)=x}, \ p_n(x(i), x)\Big).
\end{equation}
Therefore conditionally on $(x(i), {\bar I}(i),{\bar U}(i))$ we set
\begin{equation}\label{10}
\begin{array}{ll}
I^x(i+1)=&I^x(i)+ {\cal N}_n^x(i)-{\bf 1}_{I(i)>0}{\bf 1}_{x(i)=x},\\ \\
U^x(i+1)=&U^x(i)- {\cal N}_n^x(i)-{\bf 1}_{I(i)=0}{\bf 1}_{x(i)=x}
\end{array}
\end{equation}
for all $x\in S$. Then given $({\bar I}(i+1),{\bar U}(i+1))$, 
choose $v_{i+1}$ {\it uniformly } among the revealed unmarked vertices, 
unless this set is empty; in the letter  case choose $v_{i+1}$ {\it
  size-biased} among the unrevealed vertices as long as $U(j)>0$,
otherwise, stop the algorithm. In other words, we
let $x(j)=x_{v_{j}}$ have the following
 distribution
\begin{equation}\label{6}
{\bf P}\left\{x(j)=x\mid {\bar I}(j),{\bar U}(j) \right\}=
\left\{
\begin {array}{ll}
\frac{I^x(j)}{I(j)}, & \mbox{ if } I(j)> 0, \\ \\ 
\frac{xU^x(j)}{\sum_{x\in S}xU^x(j)}, & \mbox{ otherwise },
\end{array}
\right.
\end{equation}
for all $x\in S$.

We shall use Markov chain (\ref{MC}) to define the sizes of the  components.
We start with the component containing vertex $v_1$. By our algorithm 
vertex $v_i$, $i\geq 2$, belongs to the same component as $v_1$
if and only if $I(i)=\sum_{x\in S}I^x(i)> 0$. Hence, the size of the component containing $v_1$ is exactly
\begin{equation}\label{7}
\min \{j\geq 2: I(j)= 0\}-1.
\end{equation}
Given the states of  Markov chain (\ref{MC}) 
let us define a process 
which gives a useful representation for $\tau_1$ as well as for the sizes 
of other components in the graph.
Set
\begin{equation}\label{9*}
\begin{array}{ll}
z^x(1)& = \ 0,\\ \\
z^x(i+1)& = \ z^x(i)-{\bf 1}_{x(i)=x}+ {\cal N}_n^x(i) , \ \ \ i\geq 1,
\end{array}
\end{equation}
and consider
\[z(i):=\sum_{x\in S}z^x(i),   \]
which by (\ref{9*}) satisfies
\begin{equation}\label{9}
\begin{array}{ll}
z(1)& = \ 0,\\ \\
z(i+1)& = z(i)-1 +\sum_{x\in S}{\cal N}_n^x(i).  
\end{array}
\end{equation}
Notice that as long as $1< i \leq  \min \{j>1: I(j)= 0\}$ we simply have
\begin{equation}\label{8}
z(i)=I(i)-1.
\end{equation}
Indeed, by (\ref{10}) and (\ref{5})  
\begin{equation}\label{12}
\begin{array}{ll}
I(1)& = \ 0,\\ 
I(i+1)& = \ I(i)+ 
\sum_{x\in S}{\cal N}_n^x(i) -{\bf 1}_{I(i)>0}, \ \ \ 1\leq i\leq n,
\end{array}
\end{equation}
which together with (\ref{9}) gives (\ref{8}). 
Then by (\ref{8}) and  (\ref{7}) the  size of the first revealed 
component (containing $v_1$) is
\begin{equation}\label{13}
\min \{i> 1: z(i)= -1\}-1.
\end{equation}
Furthermore, setting   for all $k\geq 1$
\begin{equation}\label{14}
\tau_k := \min \{i\geq 1: z(i)= -k\},
\end{equation}
it is easy to check by induction that the 
size of
the $k$-th revealed 
component,  $k\geq 2$, is given by
\begin{equation}\label{15}
\tau_{k}- \tau_{k-1}. 
\end{equation}

By the construction of the breadth-first process  Theorem \ref{C1} 
should follow (at least intuitively) from the next stated theorem, which is the main result here.
\begin{theo}\label{T1} Define for $s\geq 0$
\begin{equation}\label{23}
Z_n(s)= n^{-1/3} \ z\left(  [ n^{2/3} s]+1\right).
\end{equation}
Then $Z_n(s)\stackrel{d}{\rightarrow} {W}^a(s)$ as $n\rightarrow \infty$.
\end {theo}

Observe that process $\{z(i), i\geq 1\}$ is not Markov. However, it
converges (after scaling) to a Markov process. The reason is that we
consider the rescaled process up to time $n^{2/3}$, and within this
time we explore, roughly speaking only small amount of vertices.

It will be proved in Lemma \ref{L1} that in our algorithm of revealing components the ordering $x(i), i\geq 1,$ is size-biased. 
By Theorem \ref{T1} the limiting process is ${W}^a(s)$, which differs from
$W(s)+as-\frac{1}{2}s^2$ considered in \cite{A} only by the constant 
coefficients. Therefore 
Theorem \ref{C1} follows from Theorem \ref{T1} by the general theory of  
 the multiplicative coalescent developed by Aldous (Lemma 14 and Proposition 15 in \cite{A}).

\bigskip
\section{Proof of Theorem \ref{T1}.}

\subsection{Weak convergence.}

For any function $f:{\bf Z}\rightarrow R$ we shall  denote
\[\Delta f(i)=f(i+1)-f(i), \ \ \ \ i\in \{1,2, \ldots\}.\]
Define a martingale sequence ${\cal M}_n(k)$, $k\geq 1$, by setting ${\cal M}_n(1)=z(1)$
and 
\begin{equation}\label{M2}
\Delta  {\cal M}_n(k) = \Delta z(k) - {\bf E} \left\{\Delta z(k)\mid {\cal F}^z_k\right\},
\end{equation}
where
\[{\cal F}^z_k= \sigma\{z (j): j\leq k\}.\]
Observe, that since  each $z(j)$ is defined as 
(Borel) function of states of the Markov chain $\left(\left( {\bar I}(i),{\bar U}(i), x(i)  \right) , i\leq j\right),$ we have
\begin{equation}\label{M3}
{\cal F}^{ {\cal M}}_{k}:=\sigma\{ {\cal M}_n(j): j\leq k\} \subseteq {\cal F}^{ z}_{k} \subseteq  {\cal F}_{k}:=\sigma\{\left({\bar U}(k), {\bar I}(k), x(k)
\right) : j\leq k\}.
\end{equation}
Then  we can write 
\begin{equation}\label{M1}
z(k)= {\cal M}_n(k)  + \sum_{j=1}^{k-1} {\bf E} \left\{\Delta z(j)\mid {\cal F}^z_j
\right\}=:{\cal M}_n(k)  + {\cal D}_n(k) .
\end{equation}
Rescale to define for all $s>0$
\[
\begin{array}{ll}
{\widetilde {\cal D}}_n (s)& = {n^{-1/3}} \ {\cal D}_n (1+[n^{2/3} s]), \\ \\
{\widetilde {\cal M}}_n (s)&  = {n^{-1/3}}
 {\cal M}_n(1+[n^{2/3} s]). 
\end{array}
\]
Our aim is to show that (Lemma \ref{L3})
\begin{equation}\label{M4}
{\widetilde {\cal M}}_n (s) 
\stackrel{d}{\rightarrow}
\sqrt{{\bf E}X {\bf E}X^3} \ W(s),
\end{equation}
while (Lemma \ref{L2}) 
\begin{equation}\label{M5}
{\widetilde {\cal D}}_n (s) \stackrel{P}{\rightarrow}as -
\frac{1}{2}
\frac{{\bf E}X^3}{{\bf E}X} s^2
\end{equation}
uniformly on a finite interval 
 as $n\rightarrow \infty$. Then Theorem \ref{T1} follows by
(\ref{M4}) and (\ref{M5}).

First we shall study the summands in ''the drift term'' ${\cal D}_n(k)$ making use of the Markov chain introduced above.

From now on we assume that $n$ is so large that 
\[-1/2 < \varepsilon _n : = \frac{a}{n^{1/3}}<1/2.\]
Also, since we are dealing here with convergence in probability or weak convergence, we may 
restrict  ourselves on the event of asymptotically high probability, where
condition (\ref{A1}), i.e., 
\begin{equation}\label{A1*}
\max_{1\leq i \leq n}x_i\leq \frac{n^\frac{1}{3}}{\omega(n)}
\end{equation}
holds. (By the assumption the probability of this event converges to
one as $n\rightarrow \infty$.)

\begin{prop}\label{P7}
For all $k \geq 1$
\begin{equation}\label{M8}
\begin{array}{ll}
{\bf E} \left\{\Delta z(k)\mid {\cal F}^z_k
\right\} & = 
  -1 +{\bf E} \left\{ x(k) \frac{1+\varepsilon _n}{n}
\sum _{x\in S} xU^x(1) \left(1-\frac{z^x(k)}{U^x(1)}\right) 
\mid {\cal F}^z_k
\right\} \\ \\
& \ \ \ 
- \frac{1+\varepsilon _n}{n} 
{\bf E} \left\{
x(k) \sum_{i=1}^{k} x(i)\left(1+ {\bf 1}_{I(i)=0}  
\right)-x(k)^2\mid {\cal F}^z_k
\right\}.
\end{array}
\end{equation}
\end{prop}

\noindent
{\bf Proof.}
By the definitions (\ref{9}) and (\ref{10}) we have
\begin{equation}\label{16}
\Delta z^x(i) = -{\bf 1}_{x(i)=x}+ 
{\cal N}_n^x(i)= -\Delta^x U(i)-{\bf 1}_{x(i)=x}-{\bf 1}_{I(i)=0}{\bf 1}_{x(i)=x}.
\end{equation}
Denote 
\[ {\bf E}_k\left\{ \cdot \right\}= {\bf E} \left\{\cdot  \mid {\cal F}_k
 \right\},\]
center $\Delta  z^x(k)$ and put 
 \begin{equation}\label{18}
\begin{array}{ll}
 M^x(1)& = \ z^x(1),\\ 
\Delta M^x(k)& =  \Delta  z^x(k)- {\bf E}_k \Delta  z^x(k)
\end{array}
\end{equation}
for all $x\in S$ and $k\geq 1$, where by (\ref{10*}) and (\ref{16}) 
\begin{equation}\label{19}
{\bf E}_k \Delta  z^x(k)= -{\bf 1}_{x(k)=x}+ (U^x(k)-{\bf 1}_{I(k)=0}{\bf 1}_{x(k)=x})p_n(x(k), x).
\end{equation}
Then  we derive recursively 
\begin{equation}\label{20}
\begin{array}{ll}
z^x(k+1) & = z^x(k) +\Delta M^x(k)+ {\bf E}_k \Delta  z^x(k) \\ \\
& = z^x(k)+\Delta M^x(k) -{\bf 1}_{x(k)=x}+ (U^x(k)-{\bf 1}_{I(k)=0}{\bf 1}_{x(k)=x})p_n(x(k), x) 
 .
\end{array}
\end{equation}
Taking into account that by (\ref{16})
\begin{equation}\label{17}
 U^x(k+1) = U^x(1)-z^x(k+1)-\sum_{i=1}^k {\bf 1}_{x(k)=x}\left(1+ {\bf 1}_{I(i)=0}  \right),
\end{equation}
we obtain from (\ref{20})
\begin{equation}\label{20A}
z^x(k+1) = z^x(k)+\Delta M^x(k) -{\bf 1}_{x(k)=x}  
\end{equation}
\[ + p_n(x(k), x) \left(U^x(1)-z^x(k)-
  \sum_{i=1}^{k-1} {\bf 1}_{x(i)=x}\left(1+ {\bf 1}_{I(i)=0} \right) -{\bf 1}_{I(k)=0}{\bf 1}_{x(k)=x}\right) .
\]
This allows us to derive 
\begin{equation}\label{M6}
\begin{array}{ll}
{\bf E} \left\{\Delta z(k)\mid {\cal F}^z_k
\right\} & = {\bf E} \left\{ \sum _{x\in S}\Delta z^x(k)\mid {\cal F}^z_k
\right\}\\ \\
 & = {\bf E} \left\{\left(\Delta z(j)-{\bf E}_k\Delta z(k)\right)
\mid {\cal F}^z_k
\right\}
  -1 \\ \\
& \ \ \ +{\bf E} \left\{ \sum _{x\in S} p_n(x(k), x) \left(U^x(1)-z^x(k)\right) \mid {\cal F}^z_k
\right\} \\ \\
& \ \ \ 
- \sum_{i=1}^{k-1}{\bf E} \left\{
p_n(x(k), x(i))\left(1+ {\bf 1}_{I(i)=0}  \right)\mid {\cal F}^z_k
\right\} \\ \\
& \ \ \ 
- {\bf E} \left\{
p_n(x(k), x(k)){\bf 1}_{I(k)=0}\mid {\cal F}^z_k
\right\} .
\end{array}
\end{equation}
Note that 
\[{\bf E} \left\{\left(\Delta z(j)-{\bf E}_k\Delta z(k)\right)
\mid {\cal F}^z_k
\right\}=0
\]
due to (\ref{M3}). Therefore  (\ref{M6}) yields
\begin{equation}\label{M7}
\begin{array}{ll}
{\bf E} \left\{\Delta z(k)\mid {\cal F}^z_k
\right\} & = 
  -1 +{\bf E} \left\{ x(k) \frac{1+\varepsilon_n}{n}\sum _{x\in S} xU^x(1) \left(1-\frac{z^x(k)}{U^x(1)}\right) \mid {\cal F}^z_k
\right\} \\ \\
& \ \ \ 
- \frac{1+\varepsilon_n}{n} {\bf E} \left\{
x(k) \sum_{i=1}^{k} x(i)\left(1+ {\bf 1}_{I(i)=0}  \right)-x(k)^2\mid {\cal F}^z_k
\right\} ,
\end{array}
\end{equation}
which proves the Proposition. \hfill$\Box$

In our analysis we will use the following result. 

\begin{prop}\label{P3} Uniformly in  $j\leq sn^{2/3}$ $($for any fixed $s>0 )$
\begin{equation}\label{J39}
\frac{1}{n} \sum _{x\in S}  U^x(j) \ {\rightarrow} \   1,
\end{equation}
\begin{equation}\label{39'}
\frac{1}{n} \sum _{x\in S} x U^x(j) \ {\rightarrow} \   {\bf E}X,
\end{equation}
and
\begin{equation}\label{40}
\frac{1}{n} \sum _{x\in S} x^2 U^x(j) \ {\rightarrow} \   {\bf E}X^2.
\end{equation}
in $L^1$  as $n\rightarrow \infty $.
\end{prop}

\noindent{ \bf Proof.} We shall establish (\ref{39'}); the rest follows exactly by the same argument.
By (\ref{10}) 
\[
U^x(j)=U^x(1)- \sum_{i=1}^{j-1}\left( 
{\cal N}_n^x(i)+{\bf 1}_{I(i)=0}{\bf 1}_{x(i)=x}
\right)
\]
for all $x\in S$. Hence
\begin{equation}\label{Ma1}
\frac{1}{n} \sum _{x\in S} x U^x(j)=\frac{1}{n} \sum _{x\in S} x U^x(1)- 
\frac{1}{n}\sum_{i=1}^{j-1} \sum _{x\in S} x\left( 
{\cal N}_n^x(i)+{\bf 1}_{I(i)=0}{\bf 1}_{x(i)=x}
\right)
\end{equation}
\[= \frac{1}{n} \sum _{x\in S} x U^x(1)- 
\frac{1}{n}\sum_{i=1}^{j-1} \sum _{x\in S} x
{\cal N}_n^x(i)
- \frac{1}{n}\sum_{i=1}^{j-1} x(i){\bf 1}_{I(i)=0},
\]
where the last term for all $j=O\left( n^{2/3}\right)$
\begin{equation}\label{Ma2}
\frac{1}{n}\sum_{i=1}^{j-1} x(i){\bf 1}_{I(i)=0}\leq 
\frac{1}{n}\sum_{i=1}^{j-1} x(i)\leq \frac{O\left( n^{2/3}\right)}{n}
\max_{1\leq i \leq n }x_i \ {\rightarrow} \ 0
\end{equation}
in $L^1$  as $n\rightarrow \infty $, since  by the assumption  (\ref{A1}) we have 
\begin{equation}\label{E}
{\bf E} \max_{1\leq i \leq n}x_i =o({n^\frac{1}{3}}).
\end{equation}
We will show that  the first summand on the right of (\ref{Ma1}) gives the main contribution.
Consider 
\begin{equation}\label{Ma4}
\frac{1}{n} \sum _{x\in S} x U^x(1)= \frac{1}{n} \sum _{x\in S} x \#\{i: x_i=x\}= \frac{1}{n} \sum _{x\in S} x \sum_{i=1}^{n} {\bf 1}_{\{x_i=x\}}
\end{equation} 
\[= \frac{1}{n} \sum_{i=1}^{n} \sum _{x\in S} x {\bf 1}_{\{x_i=x\}}.\]
Under assumption that $x_i$ are $i.i.d.$ with finite third moment, the ergodic theorem gives us  convergence in $L^1$ and $a.s.$
\[
\frac{1}{n} \sum_{i=1}^{n} \sum _{x\in S} x {\bf 1}_{\{x_i=x\}}
\ {\rightarrow} \   {\bf E} \sum _{x\in S} x {\bf 1}_{\{x_1=x\}}
= {\bf E} X,
\]
and thus by (\ref{Ma4}) we have convergence in $L^1$ and  $a.s.$
\begin{equation}\label{Ma5}
\frac{1}{n} \sum _{x\in S} x U^x(1)
\ {\rightarrow} \   {\bf E} \sum _{x\in S} x {\bf 1}_{\{x_1=x\}}
= {\bf E} X.
\end{equation}
Finally,  we bound with help of (\ref{10*})
\[
{\bf E} \frac{1}{n}\sum_{i=1}^{j-1} \sum _{x\in S} x
{\cal N}_n^x(i)\leq 
{\bf E} \frac{1}{n}\sum_{i=1}^{j-1} \sum _{x\in S}x p_n(x(i), x)U^x(i)\leq 
{\bf E} \frac{1+\varepsilon_n}{n}\sum_{i=1}^{j-1} x(i) \sum _{x\in S}x^2 \frac{U^x(1)}{n}.
\]
Using representation (\ref{Ma4}) we derive from here
\begin{equation}\label{Ma3}
{\bf E} 
\frac{1}{n}
\sum_{i=1}^{j-1} \sum _{x\in S} x {\cal N}_n^x(i)
\leq
\frac{O\left( n^{2/3}\right)}{n} 
\sum _{x\in S} x^2 \frac{1}{n} \sum _{m=1}^n {\bf E} \left(\max_{1\leq i \leq j} x(i) \right) 
{\bf 1}_{\{x_m=x\}}.
\end{equation} 
Separately we compute taking into account (\ref{E})
\[{\bf E} \left(\max_{1\leq i \leq j} x(i) \right) 
{\bf 1}_{\{x_m=x\}}
\leq {\bf E} {\bf 1}_{\{x_m=x\}}(x_m+{\bf E}\max_{1\leq i \leq n: \ i\neq m} x_i  )={\bf P} \left\{x_m=x\right\}(x+o(n^{1/3})).
\]
Substituting this into (\ref{Ma3}) we immediately derive
\[
{\bf E} \frac{1}{n}\sum_{i=1}^{j-1} \sum _{x\in S} x
{\cal N}_n^x(i) 
\leq
\frac{O\left( n^{2/3}\right)}{n} 
\sum _{x\in S} x^2 {\bf P} \left\{X=x\right\}(x+o(n^{1/3}))
\rightarrow  0.\]
This confirms that  
\begin{equation}\label{Ma6}
\frac{1}{n}\sum_{i=1}^{j-1} \sum _{x\in S} x
{\cal N}_n^x(i) \rightarrow  0
\end{equation}  
in $L^1$ uniformly in $j\leq sn^{2/3}$.
Convergence of the terms (\ref{Ma5}), (\ref{Ma6}) and (\ref{Ma2}), 
together with formula (\ref{Ma1})
yield statement (\ref{39'}). 

Statements (\ref{40}) and (\ref{J39}) can be proved exactly in the same way.
\hfill$\Box$

\begin{lem}\label{L1}
Let $x_n(i)=x(i)$, $1\leq i \leq k$, be the sequence of random variables defined in (\ref{MC}), and let $k \leq sn^{2/3}$ for some constant $s>0$. Then 
\begin{equation}\label{J9}
\begin{array}{ll}
{\bf P}\{ x_n(i) =y  \} &
=(1+o(1)){\bf E} \frac{x U^x(i)}{\sum _{x\in S}xU^x(i)}\\ \\
&
=(1+o(1)){\bf E} \frac{x U^x(1)}{\sum _{x\in S}xU^x(1)},
\end{array}
\end{equation}
where $o(1) \rightarrow 0$ as $n\rightarrow \infty $ uniformly in 
$1\leq i \leq k$, and
\begin{equation}\label{J9*}
{\bf P}\{ x_n(i) =y \mid \bar{U}(i)\} 
=(1+o_{L_1}(1)) \frac{x U^x(i)}{\sum _{x\in S}xU^x(i)},
\end{equation}
where $o_{L_1}(1) \rightarrow 0$ in ${L_1}$ as $n\rightarrow \infty $ uniformly in 
$1\leq i \leq k$.
\end{lem}

We postpone the proof of this lemma till the end of this section.

Lemma \ref{L1} together with Proposition \ref{P3} yield immediately the following corollary.
\begin{cor}\label{CL}
1. Uniformly in 
$1\leq i \leq sn^{2/3}$
\begin{equation}\label{W}
x_n(i) \stackrel{d}{\rightarrow} {\widetilde X}, \ \ \  \mbox{ as  }
n\rightarrow \infty
\end{equation}
 where ${\widetilde X}$ is a random variable with a distribution
\begin{equation}\label{dx}
{\bf P}\{ {\widetilde X}=y\}= \frac{y \, {\bf P}\{X=y\}}{{\bf E}X}.
\end{equation}

2. For all  $k \leq sn^{2/3}$
\begin{equation}\label{34A}
\frac{1}{k} \sum_{i=1}^{k} {\bf E}  x_n(i) \
{\rightarrow} \ \frac{{\bf E}X^2}{{\bf E}X},
\end{equation}
and
\begin{equation}\label{43A}
\frac{1}{k} \sum_{i=1}^{k}{\bf E} x_n^2(i){\rightarrow}  \frac{{\bf E}X^3}{{\bf E}X}
\end{equation}
 as $n\rightarrow \infty $.
\end{cor}

Now we can prove  statement (\ref{M5}). 

\begin{lem}\label{L2} Uniformly in $0 \leq s\leq s_0$
\begin{equation}\label{41}
{\widetilde {\cal D}}_n (s) \stackrel{P}{\rightarrow}as - \frac{1}{2}
\frac{{\bf E}X^3}{{\bf E}X} s^2
\end{equation}
as $n\rightarrow \infty$.
\end{lem}

\noindent
{\bf Proof.} By (\ref{M1}) and Proposition \ref{P7}
\begin{equation}\label{M11}
\begin{array}{ll}
{\cal D}_n(k+1)& = \sum_{j=1}^k {\bf E} \left\{\Delta z(j)\mid {\cal F}^z_j
\right\} \\ \\  \\
& =-k +\sum_{j=1}^k {\bf E} \left\{x(j) \frac{1+\varepsilon_n}{n}\sum _{x\in S} xU^x(1) 
\left( 1-\frac{z^x(j)}{U^x(1)}
\right)
  \mid {\cal F}^z_j
\right\} \\ \\
& \ \ \ 
- \frac{1+\varepsilon_n}{n} \sum_{j=2}^k{\bf E} \left\{
x(j) \sum_{i=1}^{j-1} x(i)\mid {\cal F}^z_j
\right\}
-\delta_n(k) ,
\end{array}
\end{equation}
where
\begin{equation}\label{M12}
\delta_n(k)=\frac{1+\varepsilon_n}{n} \sum_{j=1}^k{\bf E} \left\{x(j) \sum_{i=1}^{j} x(i)
 {\bf 1}_{I(i)=0}  \mid {\cal F}^z_j
\right\}.
\end{equation}

We shall use the following fact. 
\begin{prop}\label{ApP}
Uniformly in $j\leq s_0n^{2/3}$
\[{\bf E} \frac{|z^x(j)|}{U^x(1)} = o\left(\frac{x}{n^{1/3}}\right). \ \ \ \]
\end{prop}
\noindent
{\bf  Proof.} Observe, that
\begin{equation}\label{J10}
|z^x(j)|\leq \sum_{i=1}^{j-1}\left( \Delta ^xz(i)\right)^+ + \sum_{i=1}^{j-1}
\left( \Delta z^x(i)\right)^- ,
\end{equation}
where $a^+=\max\{0,\ a\}$ and $a^-=- \min\{0,\ a\}$.
By the definition 
(\ref{9*})
\begin{equation}\label{Ap2}
\sum_{i=1}^{j-1}\left( \Delta z^x(i)\right)^+ 
\leq \sum_{i=1}^{j-1}{\cal N}_n^x(i)
\end{equation}
and
\begin{equation}\label{Ap22}
\sum_{i=1}^{j-1}\left( \Delta z^x(i)\right)^- \leq \sum_{0\leq i\leq j-1:\ 
I(i)=0}{\bf 1}_{x(i+1)=x}.
\end{equation}
Recall also that by the definition  (\ref{6})
in our algorithm at step $i$
when $I(i)=0$ we choose $x(i+1)$ size-biased among the unrevealed vertices. Then (\ref{J10}) together with (\ref{Ap2}) and  (\ref{Ap22})  give us
\[{\bf E} 
\frac{|z^x(j)|}{ U^x(1)
} \leq 
\sum_{i=1}^{j-1} {\bf E}\frac{
{\cal N}_n^x(i)}{ U^x(1)}
+  {\bf E} 
\frac{1}{ U^x(1)
} 
\sum_{0\leq i\leq j-1:\ I(i)=0}
\frac{xU^x(i) }{\sum_{x\in S}xU^x(i)}
.
\]
Since $U^x(i)$ is non-increasing in $i$, the last bound yields
\[
{\bf E}\frac{|z^x(j)|}{U^x(1)} \leq (1+\varepsilon_n)\left(
\sum_{i=1}^{j-1} {\bf E}
\frac{ x(i)x}{n} 
+  {\bf E} 
\frac{x}{\sum_{x\in S}xU^x(j)} 
\# \{0\leq i\leq j-1:\ I(i)=0\}\right)
\]
\[\leq 2 \frac{x j}{n}\left( 
\frac{1}{j}\sum_{i=1}^{j-1} {\bf E}  x(i)
+  {\bf E} 
\frac{n}{\sum_{x\in S}xU^x(j)} \right). \]
Using Proposition \ref{P3} and Corollary \ref{CL}
we derive from here
\[{\bf E}\frac{|z^x(j)|}{U^x(1)} \leq c\frac{x j}{n}\]
for some positive constant $c$.
The statement of  Proposition \ref{ApP} follows now by the assumptions on $j.$
\hfill$\Box$

\bigskip

With a help of Proposition \ref{ApP} and bound  (\ref{A1}) we derive from (\ref{M11}) 
\[
\begin{array}{ll}
{\cal D}_n(k+1)
& =-k +(1+ o_P(1))
\sum_{j=1}^k {\bf E} \left\{x(j) \frac{1+\varepsilon_n}{n}\sum _{x\in S} xU^x(1) 
  \mid {\cal F}^z_j
\right\} \\ \\
& \ \ \ 
- \frac{1+\varepsilon_n}{n} \sum_{j=2}^k{\bf E} \left\{
x(j) \sum_{i=1}^{j-1} x(i)\mid {\cal F}^z_j
\right\}
-\delta_n(k).
\end{array}
\]
This together with Proposition \ref{P3} and Lemma \ref{L1} gives us
\[
{\cal D}_n(k+1)
 =-k +(1+\varepsilon_n+ o_P(1))({\bf E}X)
\sum_{j=1}^k {\bf E} \left\{x(j)
  \mid {\cal F}^z_j
\right\} \]
\[
- \frac{1+\varepsilon_n}{n} \sum_{j=2}^k{\bf E} \left\{
x(j) \sum_{i=1}^{j-1} x(i)\mid {\cal F}^z_j
\right\}
-\delta_n(k).
\]
Again using Lemma \ref{L1} we derive from here
\[
{\cal D}_n(k+1) =-k +{\bf E}X \left(
 k \frac{{\bf E}X^2}{{\bf E}X}+\frac{k^2}{2n}\left( \left( \frac{{\bf
         E}X^2}{{\bf E}X} \right)^3-
   \frac{{\bf E}X^2{\bf E}X^3}{({\bf E}X)^2}\right)
  \right) (1+\varepsilon_n +o_P(1))
\]
\begin{equation}\label{J14*}  
- \frac{k^2}{2n} \left( \frac{{\bf E}X^2}{{\bf E}X}\right)^2(1+o_P(1)) 
-\delta_n(k).
\end{equation}

\begin{rem}\label{cE}
Formula (\ref{J14*}) 
shows that whenever ${\bf E}X^2 \neq 1$ 
the principal term in the 
drift ${\cal D}_n(k)$ is linear in $k$, which makes out further analysis inapplicable. 
Thus we get here once again a confirmation that
the case ${\bf E}X^2=1$ is critical.
\end{rem}

Using now assumptions $\varepsilon_n =a/n^{1/3}$ and ${\bf E}X^2=1$ we get from (\ref{J14*}) 
\begin{equation}\label{M14}
{\cal D}_n(k+1) =a\frac{k}{n^{1/3}} (1+o_p(1)) 
- \frac{k^2}{2n} \frac{{\bf E}X^3}{{\bf E}X}(1+o_p(1))
-\delta_n(k) 
\end{equation}
\[=n^{1/3}\left(a\frac{k}{n^{2/3}} - \frac{k^2}{2n^{4/3}} 
    \frac{{\bf E}X^3}{{\bf E}X} \right)(1+o_p(1))
-\delta_n(k) \]
Finally, we shall find an upper bound for  $|\delta_n(k)|$ defined in (\ref{M12}). 
First we derive
\begin{equation}\label{M12A}
\delta_n(k)
 \leq 
\frac{2}{n} \sum_{j=1}^k{\bf E} \left\{x(j)\left( \max_{1\leq i \leq n}x_i
\right)
\left( \sum_{i=1}^{j} {\bf 1}_{I(i)=0} \right) \mid {\cal F}^z_j
\right\}.
\end{equation}
Recall that by our construction of process $z(i)$ 
\[\sum_{i=1}^{j} {\bf 1}_{I(i)=0} \leq - \min_{0\leq i\leq j} z(i).\]
Therefore 
 we derive from (\ref{M12A}) 
\[\delta_n (k) \leq  
 \frac{\max_{1\leq i \leq n}x_i}{n} \sum_{j=1}^k \left( \max_{0\leq i\leq j} |z(i)|\right)
  {\bf E} \left\{x(j) 
\mid {\cal F}^z_j
\right\}.
\]
This together with Lemma \ref{L1} yields for $k=O(n^{2/3})$
\begin{equation}\label{45}
\delta_n (k) \leq  (c+o_P(1))
\frac{\max_{1\leq i \leq n}x_i}{n^{1/3}} \ \max_{1\leq i\leq k} |z(i)| 
\end{equation}
for some positive constant $c$.

\begin{prop}\label{P4} 
\begin{equation}\label{46}
n^{-1/3} \max_{k \leq n^{2/3} s_0} |z(k)|  \ \mbox{ is stochastically bounded as}  \   n\rightarrow \infty.
\end{equation}
\end{prop}

\noindent
{\bf Proof.} We shall use the idea from Aldous \cite{A}.
Fix a constant $K$ and define
\begin{equation}\label{M20}
\begin{array}{ll}
T^*_n=\min\{s: |z(s)|>Kn^{1/3}\},\\ \\
T_n=\min( T^*_n, [n^{2/3} s_0]).
\end{array}
\end{equation}
Then using representation (\ref{M1})
\[z(k)= {\cal M}_n(k)  + {\cal D}_n(k),\]
we get
\begin{equation}\label{M21}
 {\bf P} \left\{ \sup_{t \leq n^{2/3} s_0} |z(t)|>Kn^{1/3}
\right\}={\bf P} \left\{  |z(T_n )|>Kn^{1/3}
\right\}
\end{equation}
\[ \leq {\bf P} \left\{  |{\cal M}_n(T_n )|+|{\cal D}_n(T_n )| >Kn^{1/3}
\right\}\]

\[ \leq \frac{{\bf E}|{\cal M}_n(T_n )|+{\bf E}|{\cal D}_n(T_n )|
 }{Kn^{1/3}}
\leq 
\frac{
\left( {\bf E}{\cal M}^2_n(T_n ) \right)^{1/2}+{\bf E}|{\cal D}_n(T_n )|
 }{Kn^{1/3}}.
\]
Next we shall find bounds for ${\bf E}{\cal M}^2_n(T_n )$ and ${\bf E}|{\cal D}_n(T_n )|$ separately. 

Consider first ${\bf E}{\cal M}^2_n(T_n )$.
Set
\begin{equation}\label{defA}
A_n (k) = \sum _{j=1}^{k} 
 {\bf E}\left\{ \left( \Delta {\cal M}_n(j) \right)^2 \mid {\cal F}^{\cal M}_{j}\right\}. 
\end{equation}
Then 
\[{\cal M}^2_n(k)-A_n (k), \ \ k\geq 1,\]
is a martingale and by the optional sampling theorem we have
\begin{equation}\label{Ap15}
{\bf E}{\cal M}^2_n(T_n ) = {\bf E} A_n (T_n) \leq  \sum _{j=1}^{[n^{2/3}s_0]} 
 {\bf E}\left( \Delta {\cal M}_n(j) \right)^2 . 
\end{equation}
By the definition 
\begin{equation}\label{Ap18}
\Delta {\cal M}_n(j)=\Delta z(j)-
 {\bf E}\left\{  \Delta z(j)\mid {\cal F}^{z}_{j}\right\},
\end{equation}
and thus 
\begin{equation}\label{Ap16}
{\bf E}\left( \Delta {\cal M}_n(j) \right)^2 = 
{\bf E}\ {\bf Var}\left\{ 
\Delta {z}(j) \mid {\cal F}^{z}_{j}\right\}.
\end{equation}
Using 
${\cal F}_{j}$ associated with the Markov chain (\ref{MC}) we derive
\begin{equation}\label{54}
{\bf Var}\left\{  \Delta z(j) \mid {\cal F}^{ z}_{j}\right\}
\end{equation}
\[
={\bf E}\left\{  
{\bf Var}\left\{ 
\Delta {z}(j) 
\mid {\cal F}_{j}\right\}
\mid {\cal F}^{z}_{j} \right\}
+{\bf Var}\left\{  
{\bf E}\left\{
\Delta {z}(j) 
\mid {\cal F}_{j}\right\}
\mid {\cal F}^{z}_{j}\right\}.
\]
Notice  
${\bf E}\left\{\Delta z(j)\mid {\cal F}^{z}_{j}\right\}$ in the definition (\ref{Ap18}) is 
${\cal F}_{j}$-measurable
 by  (\ref{M3}). Hence, we obtain from (\ref{54})
\begin{equation}\label{56}
{\bf Var}\left\{  \Delta z(j) \mid {\cal F}^{ z}_{j}\right\}
\end{equation}
\[
={\bf E}\left\{  
{\bf Var}\left\{ 
\Delta z(j) 
\mid {\cal F}_{j}\right\}
\mid {\cal F}^{z}_{j} \right\}
+{\bf Var}\left\{  
{\bf E}\left\{
\left( \Delta z(j) 
\mid {\cal F}_{j}\right\} -{\bf E}\left\{\Delta z(j)\right) 
\mid {\cal F}^{z}_{j}\right\}
\mid {\cal F}^{z}_{j}\right\}
\]
\[={\bf E}\left\{  
{\bf Var}\left\{ 
\Delta z(j) 
\mid {\cal F}_{j}\right\}
\mid {\cal F}^{z}_{j} \right\}
+{\bf Var}\left\{  
{\bf E}\left\{
\Delta z(j) 
\mid {\cal F}_{j}\right\} 
\mid {\cal F}^{z}_{j}
\right\}.
\]
Conditionally on ${\cal F}_{j}$ we have by (\ref{16})
\begin{equation}\label{57}
\Delta z(j) = \sum_{x\in S}\Delta z^x(j) =-1+ 
\sum_{x\in S}{\cal N}_n^x(j),
\end{equation}
where ${\cal N}_n^x(j)$ are independent
Binomial random variables  for different  $x$. Hence,
by (\ref{10*})
\begin{equation}\label{58}
{\bf Var}\left\{ 
\Delta z(j) 
\mid {\cal F}_{j}\right\} = (1+o(1))\frac{1}{n}\sum_{x\in S} U^x(j) x(j)x= x(j){\bf E}X(1+o_{L_1}(1)),
\end{equation}
where to derive the last equality we used Proposition \ref{P3}.
Similarly we get 
\begin{equation}\label{59}
{\bf E}\left\{
\Delta z(j) 
\mid {\cal F}_{j}\right\} = (1+o(1))\frac{1}{n}\sum_{x\in S} U^x(j) x(j)x= x(j){\bf E}X(1+o_{L_1}(1)).
\end{equation}
Substituting (\ref{58}) and (\ref{59}) into (\ref{56}) gives us
\begin{equation}\label{60}
{\bf Var}\left\{  \Delta {z}(j) \mid {\cal F}^{z}_{j}\right\}
\end{equation}
\[={\bf E}\left\{  
x(j) 
\mid {\cal F}^{z}_{j} \right\}{\bf E}X(1+o_{L_1}(1))
+{\bf Var}\left\{  
x(j) 
\mid {\cal F}^{z}_{j}
\right\} ({\bf E}X)^2(1+o_{L_1}(1)).
\]
Now with a help of (\ref{60}) we can derive from (\ref{Ap16}) for all large $n$
\[
{\bf E}\left( \Delta {\cal M}_n(j) \right)^2 \leq 2 \left(
{\bf E}  
x(j) {\bf E}X
+{\bf E}{\bf Var}\left\{
x(j) 
\mid {\cal F}^{z}_{j}
\right\} ({\bf E}X)^2\right)
\]
\[\leq 2 \left(
{\bf E}  
x(j) {\bf E}X
+{\bf E}
x^2(j) 
 ({\bf E}X)^2\right),\]
which by Lemma \ref{L1} is uniformly bounded, so that
\begin{equation}\label{Ap17}
{\bf E}\left( \Delta {\cal M}_n(j) \right)^2 \leq c
\end{equation}
for some $c>0$, which  together with (\ref{Ap15}) implies 
\begin{equation}\label{Ap19}
{\bf E}{\cal M}^2_n(T_n ) \leq  cn^{2/3}s_0 . 
\end{equation}

Let us bound now term ${\bf E}|{\cal D}_n(T_n )|$ in (\ref{M21}).
Using formula (\ref{M14}), we derive
\[
|{\cal D}_n(T_n)|\leq 
 2 |a|\frac{n^{2/3}s_0 }{n^{1/3}} 
+2 \frac{(n^{2/3}s_0 )^2}{2n} \left( \frac{1}{{\bf E}X}\right)^2
+|\delta_n(T_n-1)|, 
\]
where by (\ref{45}) and definition of $T_n$
\[
\delta_n (T_n-1) \leq  
\frac{1}{\omega_1(n)} \ \max_{1\leq i\leq T_n-1} |z(i)| \leq  
\frac{Kn^{1/3}}{\omega_1(n)} .
\]
Hence, 
\begin{equation}\label{Ap20}
{\bf E}|{\cal D}_n(T_n)| \leq c_1 n^{1/3}
\end{equation}
for some positive $c_1$.
Substituting bounds (\ref{Ap20}) and (\ref{Ap19}) into (\ref{M21}), we get
\[
 {\bf P} \left\{ \sup_{t \leq n^{2/3} s_0} |z(t)|>Kn^{1/3}
\right\}\leq 
\frac{\sqrt{cs_0}+c_1}{K},
\]
and therefore establish
(\ref{46}). \hfill$\Box$
\bigskip

Proposition \ref{P4} and bound (\ref{45}) yield
\[
\frac{\delta_n (s_0 n^{2/3})}{n^{1/3}} 
\rightarrow _P 0.\]
This together with (\ref{M14})  implies 
\[
{\widetilde {\cal D}}_n (s) = {n^{-1/3}} \ {\cal D}_n ([n^{2/3} s]+1)
\stackrel{P}{\rightarrow} as -  \frac{1}{2}\frac{{\bf E}X^3}{{\bf E}X}s^2,
\]
which completes the proof of Lemma \ref{L2}. \hfill$\Box$
\bigskip

Next we shall study the martingale part  ${\cal M}_n (k)$  in the representation  (\ref{M1}) and prove (\ref{M4}). 

\begin{lem}\label{L3} 
\begin{equation}\label{50}
{\widetilde {\cal M}}_n (s) \stackrel{d}{\rightarrow}
\sqrt{{\bf E}X {\bf E}X^3} \ W(s)
\end{equation}
as $n\rightarrow \infty$.
\end{lem}

\noindent
{\bf Proof.} To establish convergence (\ref{50}) we shall apply the functional CLT for martingales (see \cite{E}, Theorem 1.4 (b) and  Remark 1.5). 
Rescale $A_n $ defined in (\ref{defA}):
\[{\widetilde { A}}_n (s) = {n^{-2/3}}
A_n ([n^{2/3} s])=
{n^{-2/3}} \sum _{j=1}^{[n^{2/3} s]} 
 {\bf E}\left\{ \left( \Delta {\cal M}_n(j) \right)^2 \mid {\cal F}^{\cal M}_{j}\right\}, \]
where
\[\Delta{\cal M}_n (j)= \left(\Delta z(j)-
 {\bf E}\left\{  \Delta z(j)\mid {\cal F}^{z}_{j}\right\}
\right), \ \ \ \ j>1,\]
with 
\[{\cal M}_n (1)=z(1)=0.\]

We have to verify the following conditions: for each $s>0$
\begin{equation}\label{51}
\lim_{n\rightarrow \infty}{\bf E}\left[\sup_{t\leq s}|{\widetilde { A}}_n (t)- {\widetilde { A}}_n (t-)| \right]=0,
\end{equation}
\begin{equation}\label{53}
\lim_{n\rightarrow \infty}{\bf E}\left[\sup_{t\leq s}| {\widetilde {\cal M}}_n (t)- {\widetilde {\cal M}}_n (t-)|^2 \right]=0,
\end{equation}
and
\begin{equation}\label{52}
{\widetilde { A}}_n(s)  \stackrel{P}{\rightarrow}
{\bf E}X {\bf E}X^3 \ s.
\end{equation}

We start with (\ref{52}). Note that by (\ref{M3})
\[
{\widetilde { A}}_n (s) = {n^{-2/3}} \sum _{j=1}^{[n^{2/3} s]} 
{\bf E}\left\{ \left( \Delta {\cal M}_n(j) \right)^2 \mid {\cal F}^{\cal M}_{j}\right\}\]
\begin{equation}\label{61}
={n^{-2/3}} \sum _{j=1}^{[n^{2/3} s]} {\bf E}\left\{ {\bf Var}\left\{  \Delta z(j) \mid {\cal F}^{ z}_{j}\right\}\mid {\cal F}^{\cal M}_{j}\right\}.
\end{equation}
Substituting here formula 
(\ref{60}), we obtain
\begin{equation}\label{62}
{\widetilde { A}}_n (s) = {n^{-2/3}} (1+o_P(1)) \sum _{j=1}^{[n^{2/3} s]} 
\left(
{\bf E}\left\{  
x(j) 
\mid {\cal F}^{\cal M}_{j} \right\}{\bf E}X \right.
\end{equation}
\[ \left. +{\bf E}\left\{ {\bf Var}\left\{  
x(j) 
\mid {\cal F}^{z}_{j}
\right\}\mid {\cal F}^{\cal M}_{j} \right\}{\bf E}X^2
\right)
\]
\[= {n^{-2/3}} \sum _{j=1}^{[n^{2/3} s]} 
{\bf E}\left\{  
x(j) 
\mid {\cal F}^{\cal M}_{j} \right\}{\bf E}X(1+o_P(1))
\]
\[+{n^{-2/3}} \sum _{j=1}^{[n^{2/3} s]} 
\left( {\bf E}\left\{  
x(j)^2
\mid {\cal F}^{\cal M}_{j}
\right\}
-{\bf E}\left\{  
\left({\bf E}\left\{  
x(j)
\mid {\cal F}^{z}_{j}
\right\}\right)^2\mid {\cal F}^{\cal M}_{j}
\right\}
 \right)
({\bf E}X)^2(1+o_P(1))
.
\]
Using  Lemma \ref{L1}, one can derive from here
\[
{\widetilde { A}}_n (s)  =s (1+o_P(1))
\left({\bf E}X^2 + ({\bf E}X^3){\bf E}X -({\bf E}X^2)^2 \right).
\]
Under assumption (\ref{1}) that ${\bf E}X^2 =1$ this
yields 
\[{\widetilde { A}}_n (s) \rightarrow _P s{\bf E}X^3{\bf E}X ,\]
which is the statement (\ref{52}).  

Next we prove (\ref{53}). Consider
\begin{equation}\label{65}
{\bf E}
\left[\sup_{t\leq s}| 
{\widetilde {\cal M}}_n (t)- {\widetilde {\cal M}}_n (t-)|^2 
\right]
=\frac{1}{n^{2/3}}\ {\bf E}\sup_{j\leq sn^{2/3}}( \Delta {\cal M}_n(j))^2,
\end{equation}
where by the definition
\[\Delta {\cal M}_n(j)=\Delta z(j) -
{\bf E}\left\{\Delta z(j)\mid {\cal F}^{\cal M}_{j}
\right\}=
\sum_{x\in S}{\cal N}_n^x(j)-
{\bf E}\left\{\sum_{x\in S}{\cal N}_n^x(j)\mid {\cal F}^{\cal M}_{j}
\right\}.\]
Then statement (\ref{53}) is equivalent to 
\begin{equation}\label{68}\lim_{n\rightarrow \infty}
\frac{1}{n^{2/3}}{\bf E}\left(\sup_{j\leq sn^{2/3}}\left| \sum_{x\in S}{\cal N}_n^x(j)-
{\bf E}\left\{\sum_{x\in S}{\cal N}_n^x(j)\mid {\cal F}^{\cal M}_{j}
\right\}\right|\right)^2=0.
\end{equation}
Observe that $\sum_{x\in S}{\cal N}_n^x(j)$ conditionally on $x(j)$
for $j\leq sn^{2/3}$ can be well approximated by Poisson 
$Po \, (x(j){\bf E}X)$ random variable, where by Corollary \ref{CL}
$x(j) \rightarrow _d {\widetilde X}$. Therefore
we shall use the following result, which is straightforward to obtain.

\begin{prop}
Let $\xi$ be a random variable such that conditionally on ${\widetilde X}=x$
the distribution of $\xi$ is $Po \, (x{\bf E}X)$. Let $\xi_i$, $1\leq i\leq k,$ be independent copies of $\xi$, and define
\[ Y_k= \max_{1\leq i\leq k} |\xi_i- {\bf E} \xi_i|.\]
Then under condition (\ref{A1}) (i.e., when ${\bf E}{\widetilde X}^2<\infty$) we have
\begin{equation}\label{67}
\lim_{k\rightarrow \infty}
\frac{1}{k}{\bf E}Y_k^2 =0.
\end{equation}
\end{prop}
\hfill$\Box$

As a corollary to this Proposition we get convergence (\ref{68}), and therefore we establish (\ref{53}).

Finally, let us check condition  (\ref{51}). By (\ref{62})
\begin{equation}\label{69}
{\bf E}\left[\sup_{t\leq s}|A_n (t)- A_n (t-)| \right]
= {n^{-2/3}}{\bf E}
 \max_{1\leq j \leq s n^{2/3}}
\Big( {\bf E}\left\{  
x(j) 
\mid {\cal F}^{\cal M}_{j} \right\}{\bf E}X 
\end{equation}
\[+ 
\left. {\bf E}\left\{  
x(j)^2
\mid {\cal F}^{\cal M}_{j}
\right\}
-{\bf E}\left\{  
\left({\bf E}\left\{  
x(j)
\mid {\cal F}^{z}_{j}
\right\}\right)^2\mid {\cal F}^{\cal M}_{j}
\right\}
 \right)
({\bf E}X)^2(1+o_P(1))
.
\]
Using statement (\ref{J9*}) from Lemma \ref{L1} and Corollary \ref{CL}
we see that the expectation of the maximum on the right is bounded
uniformly in $n$. Therefore 
condition (\ref{51}) follows. This finishes the proof of Lemma \ref{L3}. \hfill$\Box$

\subsection{Proof of Lemma \ref{L1}.} \label{Ap}
Let us consider
$ {\bf P}\{x(i)=x\}.$ When $i=1$ by the definition (\ref{defp}) we have
\begin{equation}\label{Ma48}
{\bf P}\{x(1)=x\}={\bf E} \frac{x U^x(1)}{\sum _{x\in S}xU^x(1)},
\end{equation}
and for all $i> 1$ by  the definition (\ref{6})  we have
\begin{equation}\label{Ma12}
 {\bf P}\{x(i)=x \}=
{\bf E} \left( \frac{I^x(i)}{I(i)}{\bf 1}\{I(i)> 0\}+
\frac{xU^x(i)}{\sum_{x\in S}xU^x(i)}{\bf 1}\{I(i)= 0\}\right).
\end{equation}

\bigskip

\noindent
{\it Claim.} For all $1<i\leq sn^{2/3}$
\begin{equation}\label{Ma47}
{\bf E} \frac{I^x(i)}{I(i)}{\bf 1}\{I(i)> 0\} = (1+o(1))\ {\bf E}
\frac{xU^x(i)}{\sum_{x\in S}xU^x(i)}{\bf 1}\{I(i)> 0\}
\end{equation}
where $o(1){\rightarrow} 0$ uniformly in $i\leq sn^{2/3}$ as $n{\rightarrow} \ \infty$.

Before we prove this, let us observe that
(\ref{Ma47}) 
 together with (\ref{Ma12}) and (\ref{Ma48}) would give us 
\begin{equation}\label{Ma49}
 {\bf P}\{x(i)=x \}=(1+o(1)){\bf E} \frac{x U^x(i)}{\sum _{x\in S}xU^x(i)},
\end{equation}
where $o(1){\rightarrow} 0$ uniformly in $i\leq sn^{2/3}$ as $n{\rightarrow} \ \infty$, which is the statement of Lemma \ref{L1}.

\bigskip

\noindent
{\it Proof of the Claim.}
Recall that by  the definition (\ref{10}) for all $i> 1$ we have
\begin{equation}\label{J8}
{\bf E}  \frac{I^x(i)}{I(i)}{\bf 1}\{I(i)> 0\}
\end{equation}

\[={\bf E} \frac{I^x(i-1)+ {\cal N}_n^x(i-1)-{\bf 1}_{I(i-1)>0}{\bf 1}_{x(i-1)=x}}{
I(i-1)+ \sum_{x\in S}{\cal N}_n^x(i-1)-{\bf 1}_{I(i-1)>0}}{\bf 1}\{I(i)> 0\}.
\]
Then we observe that conditionally on 
$x(i-1)$ and ${\bar U} (i-1)$ 
the distribution of ${\cal N}_n^x(i-1)$ is binomial 
\[
{\cal N}_n^x(i-1) \in Bin \Big(U^x(i-1)-{\bf 1}_{I(1)=0}{\bf 1}_{x(i-1)=x}, \ p_n(x(i-1), x)\Big).\]
Since this can be well approximated by a Poisson distribution,  we shall make use of the following fact. 
Let $Z^x(t) \in Po\ \Big( \lambda_x \Big)$, $x\in S$, $t\geq 0$,  be independent Poisson processes with parameters $\lambda_x >0$ such that $\sum_{x\in S}\lambda_x<\infty.$ Define $Z(t) := \sum_{x\in S}Z^x(t) $ to be a superposition of these processes. Then conditionally on $Z(t)$ the distribution of  $Z^x(t)$ equals the  distribution of thinned Poisson process $Z(t)$, with a probability of thinning 
\[\frac{\lambda_x}{\sum_{x\in S}\lambda_x}.\]
In particular, this implies 
\begin{equation}\label{55}
 {\bf E} \{ Z^x(1)\mid Z(1) \}= \frac{\lambda_x}{\sum_{x\in S}\lambda_x} \ Z(1),  
\end{equation}
 which yields
\begin{equation}\label{Ma14}
{\bf E} \  \frac{Z^x(1)}{\sum_{x\in S}Z^x(1) }\ {\bf 1}\left\{ 
\sum_{x\in S} Z^x(1)> 0 \ \right\}
={\bf E} \  {\bf E} \  \left\{ \frac{Z^x(1)}{Z(1) }\ {\bf 1}\left\{ 
Z(1)> 0 \ \right\} \mid Z(1)\right\}
\end{equation}
\[= \frac{\lambda_x}{\sum_{x\in S}\lambda_x}  {\bf P}\left\{ 
Z(1)> 0 \ \right\}, \ \]
and also  for any $a, b\geq 0$
\begin{equation}\label{Ma53}
{\bf E} \  \frac{a+Z^x(1)}{b+Z(1) }\ 
{\bf 1}\left\{ 
b+Z(1)> 0 \ \right\} 
\end{equation}
\[
= {\bf E} \  \frac{a+(\lambda_x/\sum_{x\in S}\lambda_x)
Z(1)}{a+Z(1) } {\bf 1}\left\{ 
b+Z(1)> 0 \ \right\} 
. 
\]

Return to our binomial variables:  with a 
help of Poisson approximation we will show that a statement similar to (\ref{Ma53}) 
holds (up to a small error term) as well for the binomial variables ${\cal N}_n^x(i)$.

\begin{prop}\label{JP}
For any $a, b\geq 0$ 
\begin{equation}\label{J53}
{\bf E} \  \frac{a+{\cal N}_n^x(i) }{b+ \sum_{x\in S}{\cal N}_n^x(i)}\ 
{\bf 1} \left\{ 
b+\sum_{x\in S}{\cal N}_n^x(i)> 0 \ \right\} 
\end{equation}
\[
=(1+o(1)) {\bf E} \  \frac{a+ \lambda_i(x) \sum_{x\in S}{\cal N}_n^x(i) }{b+ \sum_{x\in S}{\cal N}_n^x(i)}\ 
{\bf 1} \left\{ 
b+\sum_{x\in S}{\cal N}_n^x(i)> 0 \ \right\} ,
\]
where 
\[\lambda_i(x):= \frac{x U^x(i)}{ \sum_{x\in S} xU^x(i)},\]
and $o(1) \rightarrow 0$ uniformly in $i\leq sn^{2/3}$.

In particular, if $a=b=0$
\[
{\bf E} \  
\frac{ {\cal N}_n^x(i) }{ \sum_{x\in S}{\cal N}_n^x(i)}
\ {\bf 1}\left\{ 
\sum_{x\in S} {\cal N}_n^x(i)> 0 \ 
\right\}
\]
\begin{equation}\label{J14}
=(1+o(1)) {\bf E} \ 
\frac{x U^x(i)}{ \sum_{x\in S} xU^x(i)}
 \ {\bf 1} \left\{ 
\sum_{x\in S} {\cal N}_n^x(i)
> 0 \right\},
\end{equation}
where $o(1) \rightarrow 0$ uniformly in $i\leq sn^{2/3}$.
\end{prop}

\bigskip

Before we proceed with the proof of this proposition, we shall derive 
with its help statement (\ref{Ma47}) of the Claim. In fact,
statement (\ref{Ma47}) is a particular case (when $a= b= 0$) of 
the following corollary.

\begin{cor}\label{Jc}
For any $a, b\geq 0$  we have
\begin{equation}\label{J5}
{\bf E} \  \frac{a+I^x(i) }{b+ I(i)}\ 
{\bf 1} \left\{ 
b+I(i)> 0 \ \right\} 
=(1+o(1))\  {\bf E} \  \frac{a+ \lambda_i(x) I(i) }{b+ I(i)}\ 
{\bf 1} \left\{ 
b+I(i)> 0 \ \right\} ,
\end{equation}
where 
 $o(1) \rightarrow 0$ uniformly in $i\leq sn^{2/3}$.

\end{cor}

\noindent
{\bf Proof.} We shall use 
 induction argument.
Recall that  by the definitions (\ref{5}) and  (\ref{10})
\[I^x(2)={\cal N}_n^x(1) \ \ \mbox{ and } \ \ I(2)={\sum_{x\in S}{\cal N}_n^x(1) }.
\]
Therefore, (\ref{J5}) holds for $i=2$  simply by (\ref{J53}). 

Assume now that (\ref{J5}) holds for all $2\leq i<j< sn^{2/3}$. We shall deduce that then it holds for $i=j$ as well.

By  definition (\ref{10}) we have
\begin{equation}\label{Ma52}
{\bf E} \  \frac{a+I^x(j) }{b+ I(j)}\ 
{\bf 1} \left\{ 
b+I(j)> 0 \ \right\} 
\end{equation}
\[=
{\bf E} 
\frac{a+I^x(j-1)+ {\cal N}_n^x(j-1)-{\bf 1}_{I(j-1)>0}{\bf 1}_{x(j-1)=x}}{
b+I(j-1)+ \sum_{x\in S}{\cal N}_n^x(j-1)-{\bf 1}_{I(j-1)>0}}{\bf 1}\{b+I(j)> 0\}
\]
which by (\ref{J53}) 
\[
= (1+o(1)) {\bf E} 
\frac{a+I^x(j-1)+ \lambda _{j-1}(x) \sum_{x\in S}{\cal N}_n^x(j-1)-{\bf 1}_{I(j-1)>0}{\bf 1}_{x(j-1)=x}}{
b+I(j-1)+ \sum_{x\in S}{\cal N}_n^x(j-1)-{\bf 1}_{I(j-1)>0}}{\bf 1}\{b+I(j)> 0\}.
\]
Now using (\ref{6}) we derive first
\[
{\bf E} {\bf E} \left\{
\frac{{\bf 1}_{I(j-1)>0}{\bf 1}_{x(j-1)=x}}{
b+I(j-1)+ \sum_{x\in S}{\cal N}_n^x(j-1)-{\bf 1}_{I(j-1)>0}}{\bf 1}
\{b+I(j)> 0\}\mid I(j-1), I^x(j-1)\right\}
\]
\begin{equation}\label{J6}
=
{\bf E} 
\frac{{\bf 1}
\{b+I(j)> 0\}
}{
b+I(j-1)+ \sum_{x\in S}{\cal N}_n^x(j-1)-{\bf 1}_{I(j-1)>0}}
{\bf 1}_{I(j-1)>0} \ \frac{I^x(j-1)}{ I(j-1)}.
\end{equation}
Substituting (\ref{J6}) into (\ref{Ma52}), and then using
the assumption of the induction, i.e.,  formula (\ref{J5}) for $i=j-1$, we derive
\[{\bf E} \  \frac{a+I^x(j) }{b+ I(j)}\ 
{\bf 1} \left\{ 
b+I(j)> 0 \ \right\} \]
\[=(1+o(1)) \  {\bf E} \ 
\frac{a+I^x(j-1)\left(1 -{\bf 1}_{I(j-1)>0}/ I(j-1)\right)
+ \lambda _{j-1}(x) \sum_{x\in S} {\cal N}_n^x(j-1)}{
b+I(j-1)+ \sum_{x\in S}{\cal N}_n^x(j-1)-{\bf 1}_{I(j-1)>0}}\]

\[ \times {\bf 1}\{b+I(j)> 0\}
\]

\[=(1+o(1))^2 \  {\bf E} \ 
\frac{a+ \lambda _{j-1}(x) \left(I(j-1) -{\bf 1}_{I(j-1)>0}\right)
+ \lambda _{j-1}(x) \sum_{x\in S} {\cal N}_n^x(j-1)}{
b+I(j-1)+ \sum_{x\in S}{\cal N}_n^x(j-1)-{\bf 1}_{I(j-1)>0}}\]

\[ \times {\bf 1}\{b+I(j)> 0\}
\]

\begin{equation}\label{J7}
=(1+o(1)){\bf E} \  \frac{a+\lambda _{j-1}(x)I(j) }{b+ I(j)}\ 
{\bf 1} \left\{ 
b+I(j)> 0 \ \right\}. 
\end{equation}
Recall that 
\[
\lambda_{j-1}(x)= 
\frac{x U^x(j-1)}{ \sum_{x\in S} xU^x(j-1)}
=(1+o(1))\frac{x U^x(j)}{ \sum_{x\in S} xU^x(j)}= (1+o(1)) \lambda_{j}(x),\]
which holds 
uniformly in $j\leq sn^{2/3}$ by  Proposition \ref{P3}.
Hence, we readily get from (\ref{J7}) that
\[{\bf E} 
\  \frac{a+I^x(j) }{b+ I(j)}\ 
{\bf 1} \left\{ 
b+I(j)> 0 \ \right\} 
=(1+o(1)) {\bf E} \  \frac{a+\lambda _{j}(x)I(j) }{b+ I(j)}\ 
{\bf 1} \left\{ 
b+I(j)> 0 \ \right\}, \]
which confirms (\ref{J5}) for $i=j$. This completes the proof of the Corollary,
and therefore the statement (\ref{Ma47}) follows. $\Box$

\bigskip

\noindent
{\bf
Proof of Proposition \ref{JP}.}
We shall start with statement (\ref{J14}). Note that it is enough to establish
(\ref{J14}) for $i=1$, since the difference with the general case 
is in using $\frac{x U^x(1)}{ \sum_{x\in S} xU^x(1)}$
 instead of
$$\frac{x U^x(i)}{ \sum_{x\in S} xU^x(i)}=(1+o(1))\frac{x U^x(1)}{ \sum_{x\in S} xU^x(1)}$$
 where $o(1)$ is bounded and $o(1) \rightarrow 0$ in $L_1$ uniformly in $i\leq sn^{2/3}$
by Proposition \ref{P3}.

We shall explore the following
relation between the binomial and the Poisson  distributions. 
Let $Y_{n,p} \in Bin(n,p)$ and $Z_{\lambda} \in Po(\lambda)$.
Then
it is straightforward to derive from the formulas for the corresponding
probabilities that
for all $0<p<1$  and $0\leq k\leq n$ 
\[
{\bf P}
  \{ Y_{n,p}=k\}=\frac{n!}{k!(n-k)!}p^k(1-p)^{n-k}=
  \frac{n!}{n^k(n-k)!}
\left( (1-p)e^{\frac{p}{1-p}}\right)^n e^{-n\frac{p}{1-p}}
  \frac{\left(n\frac{p}{1-p}\right)^k}{k!}
\]
\begin{equation}\label{Ma23}
= \frac{n!}{n^k(n-k)!}
\left( (1-p)e^{\frac{p}{1-p}}\right)^n {\bf P}
  \{ Z_{n\frac{p}{1-p}}=k\}.
\end{equation}

Conditionally on $\bar{U}(1)$ and $x(1)$ introduce
independent Poisson  random variables
\begin{equation}\label{Ma25}
Z^x \in Po\ \Big((U^x(1)-{\bf 1}_{x(1)=x})\ \frac{p_n(x(1), x)}{1-p_n(x(1), x)} \Big), \ \ \ \ x\in S.
\end{equation}
With a help of (\ref{Ma23}) we derive for fixed values ${\bar U}(1), x(1)$
\begin{equation}\label{Ma24}
{\bf E}\left\{  \frac{{\cal N}_n^x(1)}{\sum_{x\in S}{\cal N}_n^x(1) }{\bf 1}\left\{\sum_{x\in S}{\cal N}_n^x(1) > 0\right\}\mid {\bar U}(1), x(1) \right\}
\end{equation}
\[=  \sum_{0\leq k_y \leq U^y, \ y\in S:\newline
 \ \sum_{y\in S} k_y>0} \ \ 
\frac{k_x}{\sum_{y\in S} k_y }
\prod_{y\in S} 
\frac{\hat{U}^y!}{(\hat{U}^y)^{k_y}(\hat{U}^y-k_y)!}
\left( (1-p^y)e^{\frac{p^y}{1-p^y}}\right)^{\hat{U}^y} {\bf P}
  \{ Z^y=k_y\},
\]
where to simplify notations
we set $p^y= p_n(x(1), y)$ and $\hat{U}^y={U}^y(1) -{\bf 1}_{x(1)=y}.$
Consider first 
\[
C(n, \bar{U}, x(1)):= \prod_{y\in S}\left( (1-p^y)e^{\frac{p^y}{1-p^y}}\right)^{\hat{U}^y}
=\exp \left\{\sum_{y\in S\setminus x(1): \hat{U}^y>0} \hat{U}^y 
\left(
\frac{p^y}{1-p^y}+ \log(1-p^y) \right) 
\right\}.
\]
Notice here for a further reference, that
\begin{equation}\label{ZN}
C(n, \bar{U}, x(1))= \frac{ {\bf P}\left\{\sum_{x\in S}{\cal N}_n^x(1) = 0 \mid {\bar U}(1), x(1)\right\}
}{{\bf P}\left\{\sum_{x\in S}Z^x = 0 \mid {\bar U}(1), x(1)\right\}}.
\end{equation}
Recall that whenever $\hat{U}^y>0$ we have
$p^y=o(n^{-1/3})$. Then
\begin{equation}\label{Ma26}
C(n, \bar{U}, x(1))
=\exp \left\{\sum_{y\in S\setminus x(1): \hat{U}^y>0} \hat{U}^y 
\left(
\frac{p^y}{1-p^y}-p^y + (p^y)^2 +o((p^y)^2) \right) 
\right\}
\end{equation}
\[=
\exp \left\{\sum_{y\in S\setminus x(1): \hat{U}^y>0} \hat{U}^y
\left( 2(p^y)^2
+o((p^y)^2)\right)\right\}\]
\[=
\exp \left\{\sum_{y\in S\setminus x(1)} \hat{U}^y
\left( 2\frac{(x(1)y)^2}{n^2}
+o((p^y)^2)\right)\right\}.\]
This together with the assumption ${\bf E}X^3<\infty$ and $\max x_i= o(n^{-1/3})$ yields
\begin{equation}\label{Ma27}
C(n, \bar{U}, x(1)) = 1+ o(1)
\end{equation}
uniformly in $\bar{U}$ and $x(1)$. In particular, this together with (\ref{ZN}) gives us
\begin{equation}\label{ZN1}
{\bf P}\left\{
\sum_{x\in S}{\cal N}_n^x(1) = 0 \mid {\bar U}(1), x(1)
\right\}
=(1+ o(1)  )
{\bf P}\left\{ \sum_{x\in S} Z^x = 0 \mid {\bar U}(1), x(1)\right\}.
\end{equation}

Now we rewrite (\ref{Ma24}) with a help of (\ref{Ma27}) as
\begin{equation}\label{Ma29}
{\bf E} \left\{  \frac{{\cal N}_n^x(1)}{\sum_{x\in S}{\cal N}_n^x(1) }{\bf 1}\left\{\sum_{x\in S}{\cal N}_n^x(1) > 0\right\}\mid {\bar U}(1), x(1) \right\}
\end{equation}
\[= (1+ o(1)  )  \sum_{0\leq k_y \leq U^y, \ y\in S:\newline
 \ \sum_{y\in S} k_y>0} \ \ 
\frac{k_x}{\sum_{y\in S} k_y }
\prod_{y\in S} 
\frac{\hat{U}^y!}{(\hat{U}^y)^{k_y}(\hat{U}^y-k_y)!}
 {\bf P}
  \{ Z^y=k_y\}.
\]
We shall split the sum  on the right  into two sums, one of which  will give the major contribution while the rest will be a small term. 
Define sets 
\begin{equation}\label{Ma28}
 {B}_1:
=\left \{ 
\bar{k} =(k_y, \  y \in S):\sum_{y\in S} k_y>0, \  \right.
\end{equation}
\[\left. k_y=0 \mbox{ if } \hat{U}^y=0, \mbox{ otherwise } 0\leq k_y \leq 
\max\left\{ 1, \ \frac{\hat{U}^y}{n^{1/2}\log n} \right\}
\right\}
\]
and 
\begin{equation}\label{Ma30}
 {B}_2:=
\left\{ 
\bar{k} =(k_y, \  y \in S):  \sum_{y \in S_n} k_y>0, \ 0 \leq k_y \leq
\hat{U}^y, \ \  y\in S
\right\}  \setminus B_1.
\end{equation}
For all 
$
\bar{k} \in {B}_1$ we have
\[
\prod_{y\in S} 
\frac{\hat{U}^y!}{(\hat{U}^y)^{k_y}(\hat{U}^y-k_y)!}
 = \prod_{y\in S: \ 
1< k_y \leq 
\frac{\hat{U}^y}{n^{1/2}\log n}
} 
 \prod_{i=0}^{k_y-1}\left(1-\frac{i}{\hat{U}^y}\right)
\]
\[
= 
1 +O\left(\sum_{y\in S: \ 
1< k_y \leq 
\frac{\hat{U}^y}{n^{1/2}\log n}
}
\frac{k_y^2}{\hat{U}^y}\right)= 
1+O\left(\sum_{y\in S}
\frac{{U}^y(1)}{n\log n}
\right)=1+o(1).
\]
Hence,
\[
{\bf E}  \left\{  \frac{{\cal N}_n^x(1)}{\sum_{x\in S}{\cal N}_n^x(1) }{\bf 1}\left\{\sum_{x\in S}{\cal N}_n^x(1) > 0\right\}\mid {\bar U}(1), x(1) \right\}
\]

\[=  (1+o(1)) \left( \sum_{\bar{k} \in {B}_1}  + \sum_{\bar{k} \in {B}_2}\right)\ \ 
\frac{k_x}{\sum_{y\in S} k_y }
\prod_{y\in S} 
\frac{\hat{U}^y!}{(\hat{U}^y)^{k_y}(\hat{U}^y-k_y)!} {\bf P}
  \{ Z^y=k_y\}
\]

\bigskip

\[=   (1+o(1))^2
\sum_{\bar{k} \in {B}_1}  \frac{k_x}{\sum_{y\in S} k_y } \prod_{y\in S} {\bf P}
  \{ Z^y=k_y\}
\]
\[+ (1+o(1))\sum_{\bar{k} \in {B}_2} \ 
\frac{k_x}{\sum_{y\in S} k_y }
\prod_{y\in S} 
\frac{\hat{U}^y!}{(\hat{U}^y)^{k_y}(\hat{U}^y-k_y)!} {\bf P}
  \{ Z^y=k_y\}
\]

\bigskip

\[= (1+o(1))  
\left(
\sum_{\bar{k}: \sum_{y\in S} k_y>0}  \frac{k_x}{\sum_{y\in S} k_y }\prod_{y\in S}  {\bf P}
  \{ Z^y=k_y\}\right.\]
\[\left. -
\sum_{\bar{k} \in \{\bar{k}: \sum_{y\in S} k_y>0\} \setminus B_1}  \frac{k_x}{\sum_{y\in S} k_y }\prod_{y\in S}  {\bf P}
  \{ Z^y=k_y\}
\right)\]
\[+  (1+o(1))\sum_{\bar{k} \in {B}_2} \ 
\frac{k_x}{\sum_{y\in S} k_y }
\prod_{y\in S} 
\frac{\hat{U}^y!}{(\hat{U}^y)^{k_y}(\hat{U}^y-k_y)!} {\bf P}
  \{ Z^y=k_y\}.
\]
Hence,
\begin{equation}\label{Ma31}
{\bf E}  \left\{  \frac{{\cal N}_n^x(1)}{\sum_{x\in S}{\cal N}_n^x(1) }{\bf 1}\left\{\sum_{x\in S}{\cal N}_n^x(1) > 0\right\}\mid {\bar U}(1), x(1) \right\}
\end{equation}

\bigskip

\[= :(1+o(1))
{\bf E} \left\{  \frac{Z^x}{\sum_{y\in S} Z^y } {\bf I}
  \{ \sum_{y\in S}Z^y>0\}\mid {\bar U}(1), x(1) \right\}+R_n(x,\bar{U}, x(1)),\]
where
\[|R_n(x,\bar{U}, x(1))|\leq (1+o(1))
\sum_{\bar{k} \in \{\bar{k}: \sum_{y\in S} k_y>0\} \setminus B_1}  \frac{k_x}{\sum_{y\in S} k_y }\prod_{y\in S}  {\bf P}
  \{ Z^y=k_y\}
\]
\[+ (1+o(1)) \sum_{\bar{k} \in {B}_2} \ 
\frac{k_x}{\sum_{y\in S} k_y }
\prod_{y\in S} 
\frac{\hat{U}^y!}{(\hat{U}^y)^{k_y}(\hat{U}^y-k_y)!} {\bf P}
  \{ Z^y=k_y\}.
\]
Define a set
\[
B_3: =\left\{ \bar{k} =(k_y, \  y \in S):   k_y >
\max\left\{ 1, \ \frac{\hat{U}^y}{n^{1/2}\log n} \right\} \  \mbox{ for some } y\in S
\right\},
\]
which  contains both sets $B_2$ and
$ \{ \bar{k}: \sum_{y\in S} k_y>0 \} \setminus B_1$. Then
\begin{equation}\label{Ma32}
|R_n(x,\bar{U}, x(1))| \leq 3 \left( 1+o(1)\right)
\sum_{\bar{k} \in {B}_3} \frac{k_x}{\sum_{y\in S} k_y }\prod_{y\in S}  {\bf P}
  \{ Z^y=k_y\}.
\end{equation}
It is straightforward to derive
\[
 \sum_{\bar{k} \in {B}_3} \ \frac{k_x}{\sum_{y\in S} k_y } \
\prod_{y\in S}  {\bf P}
  \{ Z^y=k_y\} 
\]
\begin{equation}\label{Ma33}
\leq \sum_{y\in S}  {\bf P}
  \left(\{ Z^y>\max \left\{ 1, \ \frac{\hat{U}^y}{n^{1/2}\log n} \right\}\} \cap \left\{Z^x>0\right\}
\right)
\end{equation}
\[
= {\bf P}\left\{Z^x>0\right\} \sum_{x\neq y\in S }
\sum_{ k =
\max 
\left\{ 1, \ \frac{\hat{U}^y}{n^{1/2}\log n} \right\}+1
}^{\infty}
e^{-\hat{U}^y\frac{p^y}{1-p^y}}
  \frac{\left(\hat{U}^y\frac{p^y}{1-p^y}\right)^{k}}{k!}
\]
\[+ \sum_{ k =
\max 
\left\{ 1, \ \frac{\hat{U}^x}{n^{1/2}\log n} \right\}+1
}^{\infty}
e^{-\hat{U}^x\frac{p^x}{1-p^x}}
  \frac{\left(\hat{U}^x\frac{p^x}{1-p^x}\right)^{k}}{k!}\]
\[
\leq  ({\bf E} Z^x) \sum_{y\in S}\
\hat{U}^y \left( \frac{p^y}{1-p^y}\right)^2 n^{1/2}\log n \ 
\sum_{ k =2
}^{\infty}
e^{-\hat{U}^y\frac{p^y}{1-p^y}}
  \frac{\left(\hat{U}^y\frac{p^y}{1-p^y}\right)^{k-2}}{(k-1)!}
\]
\[+ \hat{U}^x \left( \frac{p^x}{1-p^x}\right)^2 n^{1/2}\log n \ 
\sum_{ k =2
}^{\infty}
e^{-\hat{U}^x\frac{p^x}{1-p^x}}
  \frac{\left(\hat{U}^x\frac{p^x}{1-p^x}\right)^{k-2}}{(k-1)!}=o(1){\bf E} Z^x.
\]
Hence, (\ref{Ma31}) together with (\ref{Ma32}) and (\ref{Ma33}) yields
\begin{equation}\label{Ma35}
{\bf E}  \left\{  \frac{{\cal N}_n^x(1)}{\sum_{x\in S}{\cal N}_n^x(1) }{\bf 1}\left\{\sum_{x\in S}{\cal N}_n^x(1) > 0\right\}\mid {\bar U}(1), x(1) \right\}
\end{equation}

\[=(1+o(1))
{\bf E} \left\{ \frac{Z^x}{\sum_{x\in S}Z^x }{\bf 1}\{\sum_{x\in S}Z^x> 0\}\mid {\bar U}(1), x(1) \right\}+o(1){\bf E} Z^x. \]
Using 
(\ref{Ma14}) we derive from here
\begin{equation}\label{Ma37}
{\bf E}  \left\{  \frac{{\cal N}_n^x(1)}{\sum_{x\in S}{\cal N}_n^x(1) }{\bf 1}\left\{\sum_{x\in S}{\cal N}_n^x(1) > 0\right\}\mid {\bar U}(1), x(1) \right\}
\end{equation}
\[= 
\frac{x (U^x(1)-{\bf 1}_{x(1)=x}) }{ \sum_{x\in S} \frac{x
(U^x(1)-{\bf 1}_{x(1)=x})}{1-p_n(x(1), x)}} {\bf P} \{ \sum_{x\in S}Z^x > 0\mid {\bar U}(1), x(1)\} 
+ o(1)
(U^x(1)-{\bf 1}_{x(1)=x})\frac{p_n(x(1), x)}{1-p_n(x(1), x)}
\] 

\[= 
\frac{x U^x(1) }{ \sum_{x\in S} x
U^x(1)} {\bf P} \{ \sum_{x\in S}Z^x > 0\mid {\bar U}(1), x(1)\} 
+ o(1)
x(1)\frac{xU^x(1)}{n}.
\]
Taking expectation on both sides, we get
\begin{equation}\label{Ma44}
{\bf E}    \frac{{\cal N}_n^x(1)}{\sum_{x\in S}{\cal N}_n^x(1) }{\bf 1}\left\{\sum_{x\in S}{\cal N}_n^x(1) > 0\right\}
\end{equation}
\[
= {\bf E}
\frac{x U^x(1) }{ \sum_{x\in S} x
U^x(1)} {\bf P} \{ \sum_{x\in S}Z^x > 0\mid {\bar U}(1), x(1) \} +R_n(x),
\]
where we denote $R_n(x)$ the remaining term, which satisfies
\begin{equation}\label{Ma41}
|R_n(x)| \leq o(1) \ {\bf E} 
x(1)\frac{xU^x(1)}{n} \leq o(1) \  {\bf E} \frac{xU^x(1)}{n}. 
\end{equation}
Recalling also the relation (\ref{ZN1}), we obtain from
 (\ref{Ma44})
\begin{equation}\label{Ma17}
{\bf E}    \frac{{\cal N}_n^x(1)}{\sum_{x\in S}{\cal N}_n^x(1) }{\bf P}\left\{\sum_{x\in S}{\cal N}_n^x(1) > 0\right\}
\end{equation}
\[
=(1+o(1))\ {\bf E}
\frac{x U^x(1) }{ \sum_{x\in S} xU^x(1)} 
{\bf E} \{ 
\sum_{x\in S}{\cal N}_n^x(1) > 0\mid {\bar U}(1), x(1) \} 
+ o(1) \  {\bf E} \frac{xU^x(1)}{n}
\]
\[
=(1+o(1))\ {\bf E}
\frac{x U^x(1) }{ \sum_{x\in S} xU^x(1)} {\bf 1}\left\{\sum_{x\in S}{\cal N}_n^x(1) > 0\right\}.\]
This completes the proof of (\ref{J14}).

Statement (\ref{J53}) can be proved in the same straightforward (but lengthy) fashion. We shall omit it here for the sake of brevity. 
This finishes the proof of Proposition \ref{JP} and therefore Lemma \ref{L1} is proved. 
$\Box$

\bigskip

\noindent
{\bf Acknowledgments.}
The author gratefully acknowledges the hospitality of the
Mittag-Leffler Institute (Djursholm, Sweden) where this work was carried out with 
partial support by the Swedish Natural Science Research Council.
The author thanks 
Remco van der Hofstad,
Svante Janson, Tomasz Luczak, and particularly Anders Martin-L{\"o}f
for the  lively discussions at the Institute, without which
this project would have been hardly completed.

\bigskip

\noindent
{\bf Remark.}
After this work was completed the author became aware of study of a
similar problem by Bhamidi, van der Hofstad and  van Leeuwaarden (see
\cite{BHL}).

\end{document}